\documentclass[12pt,leqno]{amsart}
\usepackage{amsmath}
\usepackage{graphicx,color}

\newtheorem{theorem}{Theorem}[section]

\newtheorem{lemma}[theorem]{Lemma}
\newtheorem{proposition}[theorem]{Proposition}
\theoremstyle{definition}

\theoremstyle{remark}
\newtheorem{remark}[theorem]{Remark}

\numberwithin{equation}{section}

\def\R{\mathbb{R}} %Real R
\def\C{\mathbb{C}} %Complex C
\def\B{\mathbb{B}} %Unit Ball
\begin{document}

\title{Stability of Stein structures on the Euclidean space}

\author{Herv\'e Gaussier and Alexandre Sukhov}

\address{\begin{tabular}{lll}
Herv\'e Gaussier & & Alexandre Sukhov\\
Universit\'e Joseph Fourier & & U.S.T.L. \\
100 rue des Maths & & Cit\'e Scientifique \\
38402 Saint Martin d'H\`eres & & 59655 Villeneuve d'Ascq Cedex\\
FRANCE & & FRANCE\\
 & & \\
{\rm herve.gaussier@ujf-grenoble.fr} & &
{\rm sukhov@agat.univ-lille1.fr}
\end{tabular}
}

\subjclass[2000]{32G05, 32H02, 53C15}

\date{\number\year-\number\month-\number\day}

\begin{abstract}
We give a characterization of the space $\C^n$, proving
a global version of  Newlander-Nirenberg's theorem on the integrability of complex structures
close to the standard complex structure on $\C^n$.
\end{abstract}

\maketitle

\section{Introduction}
There exist essentially two different types of Stein structures.
The first type consists of complex varieties with boundary: a model example is
provided by strongly pseudoconvex domains.  A typical example of a Stein
manifold of  second type is given by  the complex Euclidean space $\mathbb
C^n$. In many aspects (the hyperbolicity, the automorphism groups, etc.) the
properties of these two types of Stein manifolds are totally different. Our
work is motivated by  Hamilton's theorem \cite{Ha} asserting that Stein
structures of the first type are stable under small smooth  deformations (up to the
boundary)  of the complex  structure.
The focus of this paper is to prove the same sort of stability result for the
space $\mathbb C^n$.  This question is closely related to the
well-known problem of complex analysis and geometry on
the integrability of almost complex
structures.
A first fundamental  result here  is due to Newlander - Nirenberg \cite
{NewNi}. Later on many different proofs of this result  and more general results  have been
obtained by various methods (see for instance
\cite{Ch,ep-yo,Had,HiTa,Hor2,We,Z}). The main result of the present work can
be viewed as a global version of the Newlander - Nirenberg theorem   for almost complex
structures defined on the whole space $\R^{2n}$. Our approach is inspired by  H\"ormander's proof \cite{Hor2} of the
Newlander-Nirenberg theorem. We combine  H\"ormander's $L_2$-techniques with
Lempert's and Kiremidjian's results on extendibility of
CR structures \cite{Kir,Le} and with a global version of Nijenhuis-Woolf's
theorem  on the existence of pseudoholomorphic discs for certain almost
complex structures on $\R^{2n}$ .

Let $M$ be a smooth real manifold of dimension $2n$. An almost complex
structure $J$ on $M$ is a tensor field of type (1,1) on $M$ (that is a section
of $End(TM)$ satisfying $J^2 = -I$). It is called {\it integrable} if any
point in M admits an open neighborhood $U$ and a
diffeomorphism $z: U \rightarrow \mathbb B$ between $U$ and the unit ball
of $\C^n$ such that $(z_*)(J): = dz \circ J \circ dz^{-1} = J_{st}$ where $J_{st}$
denotes the standard complex structure of $\C^n$. In other words the coordinate
$z$ is biholomorphic with respect to $J$ and $J_{st}$ and $M$ admits local
complex holomorphic coordinates near every point. The Nijenhuis tensor of $J$
is defined by
$$N(X,Y) = [X,Y] - [JX,JY] + J[X,JY] + J[JX,Y].$$
A structure $J$ is called formally integrable if $N$ vanishes at every point
of $M$. The  theorem of Newlander - Nirenberg \cite{NewNi} states
that formal integrability is equivalent to integrability.

Let $J$ be an integrable smooth ($\mathcal C^\infty$) almost complex structure  on
$\R^{2n}$. We assume that for some real positive number
$\theta  > 1$  we have :
\begin{eqnarray}
\label{size0}
\parallel  D^\alpha J(z) -  D^\alpha J_{st} (z)\parallel  \leq
\frac{\lambda}{1+\parallel z \parallel^{\vert \alpha \vert +
\theta }} \ {\rm for}\ 0 \leq |\alpha | \leq 1 \ {\rm and} \ z \in
\mathbb C^n
\end{eqnarray}
and for some positive integer $K = K(n)$ :
\begin{eqnarray}
\label{size1}
\parallel D^\alpha J(z) - D^\alpha J_{st} (z)\parallel  \leq
\frac{\lambda}{1+\parallel z \parallel^{\vert \alpha \vert}} \ {\rm for}\ 2
\leq | \alpha | \leq K(n) \ {\rm and}\ z \in \mathbb C^n
\end{eqnarray}
where $\lambda > 0$ is small enough. Here and everywhere below we use the
notation $\parallel z \parallel^2 = \sum_{j=1}^n \vert z_j \vert^2$. The norm
in the left hand-side is an arbitrary fixed norm on the space of real $(2n \times
2n)$ matrices.

\vskip 0,2cm
In this paper we prove the following result.

\begin{theorem}
\label{MainTheorem}
For every $n$ there exist a positive  integer $K(n)$ and  a positive  real number $\lambda$ such that for every  smooth {\rm
  (}$\mathcal C^\infty${\rm )} integrable almost complex
structure $J$ on $\R^{2n}$ satisfying Conditions (\ref{size0}) and (\ref{size1}) the complex  manifold
$(\R^{2n},J)$ is biholomorphic to $\C^n = (\R^{2n},J_{st})$.
\end{theorem}
As we pointed out above, one can view Theorem~\ref{MainTheorem} as a global version of the classical Newlander-Nirenberg theorem \cite{NewNi} or as an analog of Hamilton's theorem \cite{Ha} for the whole complex affine space. Similar results were obtained by I.V. Zhuravlev \cite{Z}, Ch. L. Epstein and Y.Ouyang
~\cite{ep-yo} and E.Chirka \cite{Ch} under different hypothesis and by quite different techniques. We will discuss these results in the last section.

\section{The $\overline\partial$-problem on $(\C^n,J)$}

In this section we use the standard techniques of $L^2$ estimates for the $\overline\partial_J$-operator  \cite{Hor1,Hor2}. We suppose that $J$ is a smooth ($\mathcal C^\infty$) integrable almost complex structure on $\mathbb R^{2n}$ such that
\begin{eqnarray}
\label{size}
\parallel J - J_{st} \parallel_{\mathcal C^{2}(\mathbb R^{2n})} \leq \gamma,
\end{eqnarray}
where $\gamma$ is a sufficiently small real positive number.

\vskip 0,1cm Identifying
 $\R^{2n}$ with $\C^n$ we denote by $z = (z_1,...,z_n)$ the standard complex coordinates.
 Fix a basis $\omega_j$, $j=1,\dots,n,$ of
differential forms of type (1,0) with respect to $J$ on $\R^{2n}$.
Since $J$ is close to $J_{st}$ we can choose $\omega_j$ in the
form $\omega_j = dz_j + \sum_{k=1}^n b_{jk}d\overline z_k$ where
the smooth coefficients $b_{jk}$ are small enough. Then
$\overline\partial_J u = \sum_j u_{\overline \omega_j}
\overline\omega_j$ for every $u \in \mathcal C^\infty(\mathbb R^{2n})$.
Consider the Hermitian metric $ds^2 = \sum \delta_{ij}\omega_i
\overline \omega_j$ on $\R^{2n}$ (here $\delta_{ij}$ denotes the
Kroneker symbol). This metric is compatible with $J$. Since
$\overline \partial_J\overline \omega_j$ is a $(0,2)$-form, we
have $\overline \partial_J \overline\omega_j = \sum_{1 \leq k < l
\leq n}b_{kl}^j \overline \omega_k \wedge \overline \omega_l$
where the coefficients $b_{kl}^j$ can be explicitly expressed in
term of the almost complex structure tensor $J$ and its first
order derivatives. Similarly $\overline \partial_J \omega_j$ is a
$(1,1)$-form and 
\begin{eqnarray}
\label{coefficients}
\overline \partial_J \omega_j = \sum_{1 \leq k <
l \leq n}c_{kl}^j \overline \omega_k \wedge  \omega_l 
\end{eqnarray}
The smooth
coefficients $c^i_{kj}$ here  are determined by the derivatives of $J$
up to the first order.
For a real  function $\varphi$ of class $\mathcal C^\infty$ the coefficients $\varphi_{kj}$ of the Levi form of the function $\varphi$ with respect to the structure $J$ are given by :

\begin{eqnarray}
\label{Leviform}
\varphi_{kj} = \frac{\partial^2 \varphi}{\partial \overline \omega_j \partial \omega_k} + \sum_i c^i_{jk}\frac{\partial \varphi}{\partial \omega_i} = \frac{\partial^2 \varphi}{\partial \omega_k \partial \overline \omega_j} +
\sum_i \overline{c}^i_{kj}\frac{\partial \varphi}{\partial \overline \omega_i}.
\end{eqnarray}

 Denote by
$dV_0$ the standard euclidean volume form on $\R^{2n}$
 and by $dS_0$ its restriction to the boundary of $\Omega$.
 The volume form $dV = \sum_{j = 1}^n\omega_j \wedge \overline \omega_j$ defined by $ds^2$ is given by the expression
\begin{eqnarray}
\label{volume}
dV = (1 + \Phi)dV_0
\end{eqnarray}
where the density $\Phi$ can be written explicitly in terms of
$b_{jk}$.

From now on we set
$$
\varphi(z)=\parallel z \parallel^2
$$
and we suppose that $J$ is a smooth
integrable  almost complex
structure satisfying Conditions (\ref{size0}), (\ref{size1}). In particular we
assume everywhere that $\lambda$ in  Conditions (\ref{size0}), (\ref{size1})
is small enough. 

 We point out that
under these conditions there exists a positive constant $\tau_0$ such that~:
\begin{eqnarray}
\label{separation}
\tau_0\sum_{j=1}^n  dz_j  d\overline z_j \leq
\sum_{j,k=1}^n \varphi_{kj}\omega_k  \overline \omega_j.
\end{eqnarray}
Indeed, it follows from
(\ref{Leviform}) that the functions $\varphi_{jk}$ are linear combinations
of the partial derivatives of $\varphi$ up to the second order, the
coefficients of these linear combinations  consisting of the
entries  of the matrix $J- J_{st}$  and of their first order partial
derivatives (moreover,   the coefficients of the first
order partial derivatives of $\varphi$ always contain first order partial
derivatives of the entries of $J - J_{st}$). So (\ref{separation}) is a
consequence  of (\ref{size0}), (\ref{size1}).

The
expression in the left (resp.  right) hand-side of (\ref{separation}) is just the Levi form of
$\varphi$ with respect to the structure $J_{st}$ (resp. $J$). So Condition (\ref{separation}) means that the function $\varphi$ remains strictly plurisubharmonic on
$\R^{2n}$ with respect to $J$ and its Levi form is uniformly strictly positive.

The following proposition is implicitely contained in \cite{Hor1,Hor2}.

\begin{proposition}
\label{resolution} Fix a real number $\alpha$ such that  $0 < \alpha < \tau_0$. Then for every sufficiently small positive real number $\lambda$, for every smooth integrable structure $J$ on $\mathbb R^{2n}$ satisfying condition~(\ref{size1}) and for every $g \in L^2_{0,1}(\R^{2n},\varphi)$ satisfying $$\overline\partial_J g = 0$$
there exists $u \in L^2(\R^{2n},\varphi)$ such that $$\overline\partial_J u = g$$ and
\begin{eqnarray*}
\int_{\mathbb R^{2n}} \vert u \vert^2 e^{-\varphi}dV \leq 4 \int_{\mathbb R^{2n}} \vert g \vert^2 e^{-\varphi}/(\tau_0 - \alpha)dV.
\end{eqnarray*}
\end{proposition}
\proof We keep the same notations as in \cite{Hor2}. Consider the following two maps of the $\overline\partial$-complex on $(\Omega,J)$:
\begin{eqnarray*}
T = \overline\partial_J: L^2(\mathbb R^{2n},\varphi) \rightarrow L^2_{0,1}(\mathbb R^{2n},\varphi)
\end{eqnarray*}
and
\begin{eqnarray*}
S = \overline\partial_J: L^2_{0,1}(\mathbb R^{2n},\varphi) \rightarrow L^2_{0,2}(\mathbb R^{2n},\varphi),
\end{eqnarray*}
where the space $L^2_{0,k}(\mathbb R^{2n},\varphi)$, equipped with the weighted
scalar product and norm, is defined as in \cite{Hor1} or \cite{Hor2}
pp.108-109. More precisely, let a form $f$ of type $(p,q)$ be written in the form
$$f = {\sum_{\vert I \vert = p}}'{\sum_{\vert J \vert = q}}' f_{I,J}\omega^I \wedge {\overline\omega}^J$$
where $f_{I,J}$ are skew-symmetric both in $I$ and in $J$ and
$\sum'$ means that the summation is extended only over increasing
multi-indices. Then
$$\vert f \vert^2 = \frac{1}{p!q!}\sum \vert f_{I,J} \vert^2$$
and
$$\parallel f \parallel^2_\varphi = \int_{\mathbb R^{2n}} \vert f \vert^2 e^{-\varphi}dV.$$
Fix a smooth function $\eta$ with compact support in $\R^{2n}$ and equal to $1$ in a neighborhood of the origin. Then the functions
$\eta_\nu(z) = \eta(z/\nu)$ satisfy Condition
\begin{eqnarray}\label{est}
\sum_{k=1}^n \vert \partial \eta_\nu / \partial \overline \omega_k \vert \leq const.
\end{eqnarray}
Using the estimates (\ref{est}), it follows from Lemma 5.2.1 of \cite{Hor2} that the space $D_{0,1}(\C^n)$ of smooth $(0,1)$-forms with compact support in $\C^n$ is dense in
$D_{T^*} \cap D_S$ for the graph norm $$f \mapsto \parallel f \parallel_\varphi + \parallel T^*f \parallel_\varphi + \parallel Sf \parallel_\varphi$$
Following the classical method of \cite{Hor1} or \cite{Hor2} pp.77-85,107-114, we obtain for any $f \in D_{T^*} \cap D_S$~:
\begin{eqnarray}
\label{subest1}
\int_{\mathbb R^{2n}} (\tau_0 - \alpha) \vert f \vert^2 e^{-\varphi}dV \leq 4(\parallel T^* f \parallel^2_\varphi + \parallel Sf \parallel^2_\varphi)
\end{eqnarray}
 if $\lambda$ is small enough. Inequality (\ref{subest1}) enables to conclude similarly to the proof  of Theorem 4.4.1 of~\cite{Hor2}. \qed

\section{Deformation of the complex structure on $\C^n$}

Our proof of Theorem \ref{MainTheorem} is based on isotropic dilations of
coordinates ``at infinity''.  We assume everywhere that we are in the hypothesis of this Theorem.
\vskip 0,2cm

\subsection{Dilations} Our goal is to find a solution $u_k$ of  the equation
\begin{eqnarray}
\label{3.1}
\label{dilation}
\overline\partial_J u_k = \overline \partial_J z_k
\end{eqnarray}
for $k = 1,...,n$ with suitable properties. Consider  the linear
transformations (isotropic dilations at the "point at infinity") of the form~:
$$
\begin{array}{lllll}
d_\varepsilon & : & \C^n & \rightarrow \C^n\\
              &   &   z  & \mapsto \varepsilon z =: z'.
\end{array}
$$
Set $J_\varepsilon := (d_\varepsilon)_*(J)$ and $u^\varepsilon_k(z') :=
u_k(\varepsilon^{-1} z')$. If we denote by $(z_1'\dots,z_n')$ the
coordinates of $z'$ then for every $k=1,\dots,n$~:

\begin{eqnarray}
\label{3.2}
\overline\partial_{J_\varepsilon} u_k^\varepsilon = \overline \partial_{J_\varepsilon} z'_k.
\end{eqnarray}

\vskip 0,1cm
Since $J_\varepsilon(z') = J(\varepsilon^{-1}z')$ and $J_{st}(z') = J_{st}(\varepsilon^{-1}z')$,
 Conditions (\ref{size0}) and (\ref{size1}) imply ~:

\begin{eqnarray}
\label{size20}
\parallel  D^\alpha J_\varepsilon(z') -  D^\alpha J_{st}(z') \parallel \leq
\frac{\varepsilon^\theta\lambda }{\varepsilon^{\vert \alpha \vert
+ \theta} +
\parallel z'\parallel^{\vert \alpha \vert + \theta}}, \ {\rm for \ every}\ z' \in
\mathbb C^n, 0 \leq \vert \alpha \vert \leq 1
\end{eqnarray}
and
\begin{eqnarray}
\label{size2}
\parallel D^\alpha J_\varepsilon(z') - D^\alpha J_{st}(z') \parallel \leq
\frac{\lambda }{\varepsilon^{\vert\alpha\vert} +
\parallel
  z'\parallel^{\vert \alpha \vert}} \ {\rm for \ every}\ z' \in
  \mathbb C^n, 2
  \leq  |\alpha| \leq K(n).
\end{eqnarray}
In particular for $\varepsilon \leq 1$~:
\begin{eqnarray}\label{size30}
\parallel  D^\alpha J_\varepsilon(z') - D^\alpha J_{st}(z')
\parallel \leq 2^{\vert \alpha \vert + \theta}\lambda, \ {\rm for}
\ z' \in \mathbb C^n \backslash \frac{1}{2}\mathbb B, \vert \alpha
\vert \leq 1
\end{eqnarray}
and

\begin{eqnarray}\label{size3}
\parallel D^\alpha J_\varepsilon(z') - D^\alpha J_{st}(z')
\parallel \leq
2^{\vert\alpha\vert}\lambda, \ {\rm for} \ 2 \leq |\alpha | \leq
K(n) \ \ \displaystyle z' \in \mathbb C^n \backslash
\frac{1}{2}\mathbb B.
\end{eqnarray}

Unfortunately we cannot control the behavior of the dilated structure $J_\varepsilon$ near the
origin. For this reason we modify this structure on the ball
$\frac{1}{2}\mathbb B$.

\vskip 0,1cm
\begin{remark}
Here we have to warn the reader that a given integrable CR structure
on the unit sphere may be extended to different (integrable)
smooth complex structures on the unit ball as shows the
following elementary example. Consider a real smooth
diffeomorphism $f$ of the closed unit ball in $\C^n$ which is
equal to the identity mapping in a neighborhood of the unit
sphere. Then the standard CR structure of the unit sphere, i.e.
the restriction of $J_{st}$ to the unit sphere, may be extended to
both the integrable complex structures $J_{st}$ and
$f_\star(J_{st})$ on the unit ball. Finally we may choose $f$ such that
$(f)_\star(J_{st})$ is different from $J_{st}$ in the ball.
\end{remark}

\vskip 0,1cm
In our situation, consider  the restriction   ${J_\varepsilon}_{\vert \frac{1}{2}\mathbb
S}$ of the complex structure $J_\varepsilon$ to the real sphere $\frac{1}{2}\mathbb
S$ where $\mathbb S$ denotes the real unit sphere of $\C^n$. This is a smooth ($\mathcal C^\infty$) CR structure
on the real sphere $\frac{1}{2}\mathbb S$. Furthermore  according to
Conditions~(\ref{size30}) and (\ref{size3}) this is a  $\mathcal C^{K(n)}$ deformation of
the standard CR structure of the sphere, namely of the
restriction of $J_{st}$ on $\frac{1}{2}\mathbb S$. It follows from the
well-known results of G.Kiremidjian \cite{Kir} (in the case $n \geq 3$) and  L.Lempert \cite{Le} (in the case $n = 2$) that for  $K(n)$ large
enough and for $\lambda$ small enough we can extend
the CR structure
${J_\varepsilon}|\frac{1}{2}\mathbb S$ to an integrable almost
complex structure $\tilde{J}_\varepsilon$ of class $\mathcal C^\infty$
on $\mathbb C^n$, coinciding with $J_\varepsilon$ on $\mathbb C^n
\backslash \frac{1}{2}\mathbb B$ and satisfying the condition~:

\begin{equation}\label{size3'}
\parallel D^\alpha \tilde J_\varepsilon(z') - D^\alpha J_{st}(z')
\parallel \leq C \lambda,
\end{equation}
for all $z' \in \frac{1}{2}\mathbb B$ and  $ |\alpha | \leq K(n)$
where $C$ is a real positive constant independent of $\varepsilon$.

\begin{lemma}\label{lemma-inv}
Condition (\ref{separation}) holds in the coordinates $z'$
for the function $\varphi: z' \mapsto  \parallel z'\parallel^2$ and
the structure $\tilde J_\varepsilon$ with $\tau_0$ independent of
$\varepsilon$.
\end{lemma}
\proof Denote by ${\mathcal L}^J(\psi,p,v)$ the value of the Levi form of a
function $\psi$ at $(p,v) \in T(\R^{2n})$.
Then Condition (\ref{separation}) can be written
\begin{eqnarray}
\label{separation1}
{\mathcal L}^{J_{st}}(\varphi,p,v) \leq {\mathcal L}^J(\varphi,p,v).
\end{eqnarray}
Set $\varphi_\varepsilon(z') = \varphi \circ (d_\varepsilon)^{-1}(z')$ that is
$\varphi_\varepsilon(z') = \varepsilon^{-2} \varphi(z')$. The invariance of the Levi
form and the inequality (\ref{separation1}) then imply
\begin{eqnarray*}
{\mathcal L}^{J_{st}}(\varphi,p',v') \leq {\mathcal L}^{J_\varepsilon}(\varphi,p',v')
\end{eqnarray*}
for any $(p',v') \in T(\R^{2n})$. Since the structure $\tilde J_\varepsilon$
is close to $J_{st}$ in the $\mathcal C^2$-metric on the ball
$(1/2){\mathbb B}$, the
last inequality also holds for $\tilde J_\varepsilon$. This proves Lemma \ref{lemma-inv}. \qed

\bigskip
The form $\bar{\partial}_{\tilde{J}_\varepsilon}z'_k$ satisfies, for every real vector field $Y$~:

$$
\begin{array}{lcl}\label{eqn1}
\bar{\partial}_{\tilde{J}_\varepsilon}z'_k(Y) &:=&
\frac{1}{2}[dz'_k(Y)+\tilde{J}_\varepsilon(z'_k)dz'_k(iY)]\\
 & = &
\bar{\partial}_{J_{st}}z'_k(Y)+
\frac{1}{2}[\tilde{J}_\varepsilon(z'_k)dz'_k(iY)-J_{st}(z'_k)dz'_k(iY)]\\
& = &
\frac{1}{2}[\tilde{J}_\varepsilon(z'_k)dz'_k(iY)-J_{st}(z'_k)dz'_k(iY)].
\end{array}
$$

It follows
from the construction of $\tilde{J}_\varepsilon$ (see
Condition~(\ref{size3})) that there exists a positive constant $M$
independent of $\varepsilon$ such that for every $k=1,\dots,n$ and for every $0 \leq |\alpha| \leq K(n)$~:
\begin{eqnarray}\label{ESTI}
\parallel D^\alpha(\overline\partial_{\tilde J_\varepsilon} z'_k)
\parallel_{L^\infty(\C^n)} \leq M.
\end{eqnarray}
Using (\ref{ESTI}) with $|\alpha|=0$  there exists a positive constant $M'$ such that~:

\begin{eqnarray*}
\parallel \overline\partial_{\tilde J_\varepsilon} z'_k \parallel_\varphi^2
\leq M'.
\end{eqnarray*}
Denote by $dV_\varepsilon$ the volume form
corresponding to the structure $\tilde J_\varepsilon$, given by
(\ref{volume}) (with $\Phi$ being defined by means of $\tilde
J_\varepsilon$).

Fix $1 \leq k \leq n$ and  apply Proposition
\ref{resolution} to the structure $\tilde J_\varepsilon$, with
$g=\overline\partial_{\tilde J_\varepsilon}z'_k$, by virtue of the inequalities (\ref{size3}) and (\ref{size3'}). There exists a solution
$u^k_\varepsilon$ of the equation (\ref{3.2}) satisfying~:

\begin{eqnarray}\label{ESTIM}
\int_{\C^n} \vert u_k^\varepsilon \vert^2 e^{-\varphi} dV_\varepsilon
\leq \int_{\C^n} \vert \overline\partial_{\tilde J_\varepsilon} z'_k \vert^2
e^{-\varphi}/\tau_1 dV_\varepsilon = (1/\tau_1) \parallel
\overline\partial_{\tilde J_\varepsilon} z'_k \parallel_\varphi^2,
\end{eqnarray}
with $\tau_1:= \tau_0 - \alpha$ for $\tau_0$ and $\alpha$ from Proposition
\ref{resolution} independent of $\varepsilon$.

Furthermore it follows from Conditions (\ref{size30}), (\ref{size3}) and (\ref{size3'}) and from the expression of $\bar{\partial}_{\tilde{J}_\varepsilon}z'_k$ given above that
for every compact subset $X$ of $\C^n$ we have~:
\begin{eqnarray}\label{EESSTT}
\lim_{\lambda \rightarrow 0}
\sup_{z \in X, 0 \leq |\alpha| \leq K(n)}|D^\alpha(\bar{\partial}_{\tilde{J}_\varepsilon}z'_k)|=0,
\end{eqnarray}
uniformly with respect to $\varepsilon$ sufficiently small.

Then using (\ref{EESSTT}) with $|\alpha| = 0$ we have by the Lebesgue
Theorem~:

\begin{eqnarray}\label{size10}
\lim_{\lambda \rightarrow 0}\|\bar{\partial}_{\tilde{J}_\varepsilon}z'_k
\|^2_\varphi = 0,
\end{eqnarray}
uniformly with respect to $\varepsilon$ sufficiently small.

It follows now from (\ref{ESTIM}) and (\ref{size10})~:
\begin{eqnarray*}
\lim_{\lambda \rightarrow 0}
\int_{\C^n} \vert u_k^\varepsilon \vert^2 e^{-\varphi} dV_\varepsilon = 0
\end{eqnarray*}
which implies~:
\begin{eqnarray}\label{size11}
\lim_{\lambda \rightarrow 0}
\int_{2\mathbb B} \vert u_k^\varepsilon \vert^2  dV_\varepsilon = 0,
\end{eqnarray}
uniformly with respect to $\varepsilon$ sufficiently small.

Finally according to (\ref{EESSTT}) we have for every $0 \leq |\alpha| \leq K(n)$~:
\begin{eqnarray}\label{ESTIMA}
\lim_{\lambda \rightarrow 0}\int_{2\mathbb B}
|D^{\alpha}(\bar{\partial}_{\tilde{J}_\varepsilon}u^\varepsilon_k)|^2dV_\varepsilon =
\lim_{\lambda \rightarrow 0}\int_{2\mathbb B}
|D^{\alpha}(\bar{\partial}_{\tilde{J}_\varepsilon}z'_k)|^2dV_\varepsilon=0,
\end{eqnarray}
uniformly with respect to $\varepsilon$ sufficiently small.

\subsection{Estimate of the ${\mathcal C}^1$-norm} Since the right hand-side
of the equality (\ref{dilation}) is smooth, according to the
well-known results (for instance, Theorem 5.2.5 of \cite{Hor2}) the
solution $u^\varepsilon_k$ is also smooth. In order to control its
${\mathcal C}^1$ norm we follow classical arguments
of~\cite{Hor2}. The next statement is just a consequence of the
Sobolev embedding Theorem~:

\begin{lemma}
\label{lem-ref}
There exists a positive constant $C_1$ independent of $\varepsilon$
such that
\begin{eqnarray*}
\parallel u_k^\varepsilon \parallel^2_{\mathcal C^1(\B)} \leq C_1
  \left ( \int_{2\B} \vert u^\varepsilon_k \vert^2dV_\varepsilon + \int_{2\B} \vert \overline \partial_{\tilde J_\varepsilon}
  {u}_k^\varepsilon \vert^2 dV_\varepsilon + \sum_{1 \leq \vert \alpha \vert
  \leq n+1} \int_{2\B} |D^\alpha (\overline\partial_{\tilde
  J_\varepsilon} u^\varepsilon_k)|^2dV_\varepsilon \right ),
\end{eqnarray*}
where $\B$ denotes the unit ball in $\C^n$.
\end{lemma}

\proof  Fix a smooth $\mathcal C^\infty$ function $\xi$ in $\C^n$ with
compact support in $2\B$, and such that the restriction of $\xi$ to
$\B$ is identically equal to $1$.  The classical arguments of \cite{Hor2}, Lemma 5.7.2 page 140 shows that there exists a positive constant $C_3$
independent of $\varepsilon$ such that~:
\begin{eqnarray}
\label{parts}
\parallel \partial_{\tilde J_\varepsilon} (\xi u_k^\varepsilon)
\parallel_{L^2(2\B)} \leq C_3 \left(\parallel (\xi u_k^\varepsilon)
\parallel_{L^2(2\B)}+ \parallel \overline\partial_{\tilde
J_\varepsilon}(\xi u_k^\varepsilon) \parallel_{L^2(2\B)}\right)
\end{eqnarray}
where the $L^2$-norms are considered with respect to the standard
volume form $dV_0$ on $\C^n$. 

It follows from (\ref{parts}) and from
Proposition~\ref{resolution} that there exist positive constants
$C_3,\ C_4$ and $C_5$ such that~:
$$
\begin{array}{lll}
\parallel \partial_{\tilde J_\varepsilon} (\xi u_k^\varepsilon)
\parallel_{L^2(2\B)} & \leq & C_3 \left(\parallel (\xi
u_k^\varepsilon) \parallel_{L^2(2\B)}+ \parallel u_k^\varepsilon
\parallel_{L^2(2\B)} \times \parallel \overline\partial_{\tilde
J_\varepsilon} \xi \parallel_{L^2(2\B)} + \parallel
\overline\partial_{\tilde J_\varepsilon} u_k^\varepsilon
\parallel_{L^2(2\B)}\right)\\ & \leq & C_4 \left(\parallel
u_k^\varepsilon \parallel_{L^2(2\B)}+ \parallel
\overline\partial_{\tilde J_\varepsilon} u_k^\varepsilon
\parallel_{L^2(2\B)}\right)\\ & \leq & C_5 \parallel
\overline\partial_{\tilde J_\varepsilon} u_k^\varepsilon
\parallel_{L^2(2\B)}.
\end{array}
$$

Thus we have~:
$$ \parallel \partial_{\tilde J_\varepsilon} ( u_k^\varepsilon)
\parallel_{L^2(\B)} \leq C_5 \parallel \overline\partial_{\tilde
J_\varepsilon} u_k^\varepsilon \parallel_{L^2(2\B)}.
$$

In particular, we may impose that the $L^2$-Sobolev norm
$W^{1,2}(2\B)$ of $u^k_\varepsilon$ is arbitrarily small if $\lambda$,
provided by Condition~(\ref{size1}), is small enough. In order to
obtain an estimate of the $\mathcal C^1$ norm, we iterate this
argument.  For instance, set $g^k_i = \partial
u_k^\varepsilon/\partial \omega_i$. Then the above argument gives the
existence of positive constants $C_6$ and $C_7$ such that~:
\begin{eqnarray*}
& &\parallel \partial_{\tilde J_\varepsilon} g^k_i \parallel_{L^2(\B)}
\leq C_6 \left(\parallel g^k_i \parallel_{L^2(2\B)}+ \parallel
\overline\partial_{\tilde J_\varepsilon} g^k_i
\parallel_{L^2(2\B)}\right) \leq C_7 \sum_{\vert\alpha \vert \leq
1}\parallel D^\alpha \overline\partial_{\tilde J_\varepsilon}
u_k^\varepsilon\parallel_{L^2(2\B)}.
\end{eqnarray*}
Iterating this argument we obtain estimates on the $L^2$-norms of the
derivatives of $u^k$ up to order $n+1$. Then the Sobolev embedding theorem
implies the desired statement. \qed

\section{$J$-complex curves and  Montel's property}
In this section  we  conclude the proof of Theorem~\ref{MainTheorem}.

\subsection{Montel's property}

Let $r$ be a sufficiently small real positive number.
According to Conditions (\ref{size11}), (\ref{ESTIMA}) and to Lemma~\ref{lem-ref},
if $K(n) \geq n+1$ then
we can choose $\lambda$
small enough such that~:
\begin{equation}\label{size4}
\parallel u_k^\varepsilon \parallel_{\mathcal C^1(\B)} \leq r,
\end{equation}
uniformly with respect to $\varepsilon$ sufficiently small.

Consider the smooth $\mathcal C^\infty$ almost complex structure
$\tilde J$ defined on $\C^n$ by $\tilde J = (d_\varepsilon^{-1})_*(\tilde
J_\varepsilon)$. In our situation we have
$\tilde J(z)=\tilde
J_\varepsilon(\varepsilon z)$. Then the structure $\tilde J$ coincides with $J$ on
$\C^n \backslash \frac{1}{2\varepsilon}\B$. In particular the function
$\tilde v_k^\varepsilon$ defined by $\tilde v_k^\varepsilon(z) :=
u_k^\varepsilon(\varepsilon z)$ satisfies the equation~: $$
\overline\partial_{J} \tilde v_k^\varepsilon = \overline \partial_{J}
z_k$$ on $\C^n \backslash \frac{1}{2\varepsilon}\B$
This means that
the function $$f^\varepsilon_k
:= z_k - \tilde v_k^\varepsilon$$ is  $J$-holomorphic on $\C^n \backslash
\frac{1}{2\varepsilon}\B$. Since the unit ball $\mathbb B$ equipped with the
structure $J$ is a Stein manifold, by the
removal of compact singularities the function $f^\varepsilon_k$ extends  on
$\C^n$  to a function still denoted
by $f^\varepsilon_k$ and $J$-holomorphic on $\C^n$.
Moreover
according to the estimate~(\ref{size4}) we have~:
\begin{eqnarray}\label{EST}
\parallel \tilde{v}_k^\varepsilon \parallel_{\mathcal C^1((1/2\varepsilon)\mathbb S)} \leq r.
\end{eqnarray}
Consider now a sequence $(\varepsilon_j)_j$ decreasing to $0$ and set
$f^j:=f^{\varepsilon_j}$.
It follows from (\ref{EST})~:
\begin{eqnarray}
\label{max}
\parallel f^j - Id  \parallel_{\mathcal C^1((1/2\varepsilon_j){\mathbb S}} \leq r.
\end{eqnarray}

In order to prove that the sequence $(f^j)$ is a normal family it is
sufficient to show that the sequence $(f^j)$ is uniformly bounded on
$r_0\mathbb B$ for every $r_0 > 0$. A natural
idea is to use the estimate (\ref{max}) and
the maximum principle. The obstacle here  is that the identity map $Id:
(\C^n,J) \longrightarrow (\C^n,J_{st})$ is not
holomorphic. For this reason  we use the techniques of $J$-holomorphic curves
similarly to  M.Gromov's work \cite{Gr} (basic local results were obtained
in \cite{NiWo}).

In what follows we denote by $\Delta$ the unit disc of $\C$.
 We remind that for $R > 0$ a smooth map  $L_J:R\Delta \rightarrow
\C^n$ is called $J$-holomorphic if  $dL_J \circ J_{st} = J \circ
dL_J$. In other words, $L_J$ is holomorphic with respect to the
structure $J_{st}$ on $\C$ and the structure $J$ on $\C^n$. We
also call such maps {\sl $J$-holomorphic discs}. We always suppose that
a $J$-holomorphic disc is continuous on $R\overline\Delta$.
In the case where a $J$-holomorphic map $L_J$ is defined on the whole complex plane $\mathbb C$ we call it {\sl $J$-complex line}.

\begin{proposition}
\label{discs}
There exists a positive constant $C$ such that for every $v \in \mathbb S$ there
is a $J$-complex line $L_J^{v} : \mathbb C \rightarrow \mathbb C^{n}$ satisfying~:
$$\sup_{\zeta \in \mathbb C}\parallel L_J^{v}(\zeta) - \zeta v \parallel \leq C \lambda,$$
$\lambda$ being a constant from (\ref{size0}), (\ref{size1}). Furthermore
$$
\mathbb C^n \backslash \mathbb B \subset \{L_J^v(\zeta), \ v \in \mathbb S,\ \zeta \in \mathbb C\}.
$$
\end{proposition}
A weaker version (but still sufficient for the proof of Lemma~\ref{bound}) of Proposition~\ref{discs} giving the existence of arbitrary large $J$-complex discs instead of lines follows from the results of \cite{c-g-s}. On the other hand Proposition~\ref{discs} is a consequence of a more general result, concerning
smooth almost complex deformations of the standard complex structure.
Since this result is of independent interest we state and prove it in the next section for the convenience of the reader (see Theorem~\ref{deform}). We will apply Proposition~\ref{discs} together with the following Lemma. This frequently used statement is just a variation of the classical Nijenhuis-Woolf's Theorem (see \cite{Si}).

\begin{lemma}\label{nw}
For $a \in \mathbb C^{n-1}$, $|a| \leq 1$, let the map $N_a : 2\Delta \rightarrow \mathbb C^n$ be defined by $N^a(\zeta) = (a,0) + \zeta (0,1)$ where $(a,0),\ (0,1) \in \mathbb C^{n-1} \times \mathbb C$. For every $\lambda$ sufficiently small there exists a $J$-holomorphic curve $N^a_J : 2\Delta \rightarrow \mathbb C^n$ such that $N^a_J(0) = (a,0)$ and $\sup_{\zeta \in 2\Delta}\|N^a_J(\zeta)-N^a(\zeta)\| \leq C \lambda$. These curves form a foliation of a neighborhood of $\bar{\mathbb B}$.

\end{lemma}

We continue now the proof of Theorem~\ref{MainTheorem}.

The crucial consequence of Proposition \ref{discs} is the following

\begin{lemma}
\label{bound}
We have
\begin{equation}\label{EST1}
\parallel f^j(z) - z \parallel \leq r + C\lambda
\end{equation}
 for every $j$ and every $z \in (1/2\varepsilon_j)\mathbb B$. In particular,
for every $r_0 > 0$ the sequence $(f^j)$ is uniformly bounded on the ball
$r_0\mathbb B$.
\end{lemma}
\vskip 0,2cm
\proof Let $p \in (1/2\varepsilon_j)\mathbb B \backslash \bar{\mathbb B}$. Consider the $J$-complex line $L_J^v$ given by Proposition~\ref{discs} passing through $p$, i.e. satisfying $L_J^v(\zeta_p)=p$ for some $\zeta_p \in \mathbb C$. Denote by $L^v$ the usual complex line $\zeta \in \mathbb C \mapsto \zeta v$. Since the maps
 $L_J^{v}: (\mathbb C,J_{st}) \mapsto (\C^n,J)$ and $f^j: (\C^n,J) \rightarrow
 (\C^n,J_{st})$ are holomorphic, the map $\zeta \mapsto f^j(L_J^{v}(\zeta))
- L^{v}(\zeta)$ is holomorphic with respect to $J_{st}$.
It follows by Proposition \ref{discs} that the domain $D_j :=
(L_J^{v})^{-1}((1/2\varepsilon_j)\mathbb B)$ is a bounded domain
in $\mathbb C$ containing the point $\zeta_p$. Furthermore, the function
$\varphi(z) = \parallel z
\parallel^2$ is $J$-plurisubharmonic and so the composition $\varphi \circ
L_J^{v}$ is a
 subharmonic function on $\mathbb C$. Applying the maximum principle to this function we obtain that the domain $D_j$ is simply connected and its boundary
 coincides with the level set $\{\zeta \in \mathbb C : \parallel L_J^{v}(\zeta)
 \parallel = 1/(2\varepsilon_j)\}$.
For $\zeta \in \C$ satisfying $\parallel L_J^{v}(\zeta)
 \parallel = 1/(2\varepsilon_j)$ we have~:
\begin{eqnarray*}
\parallel f^j(L_J^{v}(\zeta)) - L^{v}(\zeta) \parallel \leq \parallel f^j(L_J^{v}(\zeta)) -
L_J^{v}(\zeta) \parallel + \parallel L_J^{v}(\zeta) - L^{v}(\zeta) \parallel.
\end{eqnarray*}
The first quantity in the right hand-side is bounded from above by $r$ in view of
(\ref{max}) and the second one is bounded by $C\lambda$. Applying on the
domain $D_j$ the
maximum principle to the map $\zeta \mapsto f^j(L_J^{v}(\zeta))
- L^{v}(\zeta)$ we conclude that the estimate~(\ref{EST1}) holds on $(1/2\varepsilon_j)\mathbb B \backslash \bar{\mathbb B}$.

If $p \in \bar{\mathbb B}$ we repeat the same argument, replacing respectively $L^v$ by $N^a$ and $L^v_J$ by $N^a_J$, from Lemma~\ref{nw}. This ends the proof of Lemma~\ref{bound}. \qed

\subsection{ End of the proof of Theorem \ref{MainTheorem} } Denote by
$Jac(f^j)(z)$ the Jacobian determinant of the map $f^j$ at $z$.

\begin{lemma}\label{lemA}
There exists a positive constant $\alpha$ such that
$$\vert Jac(f^j)(z) \vert \geq \alpha$$
for every $j$ and every $z \in (1/2\varepsilon_j)\mathbb B$.
\end{lemma}

\proof

It follows from (\ref{max}) that the Jacobian determinant $Jac(f^j)$ of $f^j$ is a
holomorphic (with respect to $J$ and $J_{st}$) function not vanishing on $(1/2\varepsilon_j)\mathbb S$.
Hence it follows from the removal compact singularities theorem that
$Jac(f^j)$ does not vanish on $(1/2\varepsilon_j)\mathbb B$. Applying the
maximum principle to the function $1/Jac(f^j)$ we conclude. \qed

\vskip 0,2cm Lemma \ref{bound} implies that the sequence $(f^j)$
contains a subsequence (still denoted by $(f^j)$) uniformly
converging (with all partial derivatives) on any compact subset of
$\C^n$. Denote by $f$ the limit map.  Then  the map $f$ satisfies
the equation $\overline\partial_J f = 0$ that is $f: (\C^{n},J)
\rightarrow (\C^{n},J_{st})$ is holomorphic. Furthermore according
to Lemma~\ref{lemA}  the map $f$  is locally biholomorphic. On the
other hand, it follows by Lemma \ref{bound} that for every compact
subset $K$ of $\C^n$ we have $\parallel f - Id
\parallel_{\mathcal C^o(K)} \leq r + C\lambda$ and therefore
$f:\C^n \longrightarrow \C^n$ is a proper map. Now by the
classical theorem of J.Hadamard \cite{Had} $f$ is a global
diffeomorphism of $\R^{2n}$ and so is globally biholomorphic.

This completes the proof of
Theorem \ref{MainTheorem}. \qed

\section{Stability of curves under a deformation of $J_{st}$}
In this Section we establish a general version of
Proposition~\ref{discs}. This does not use the integrability of
$J$ and is valid for any almost complex structure on $\C^n$
satisfying Conditions (\ref{size0}), (\ref{size1}). We point out
that one of the results of Gromov \cite{Gr}  gives the existence
of a $J$-complex map from the Riemann sphere to  the complex
projective space $\C \mathbb P^n$ equipped with an almost complex
structure tamed by the standard symplectic form of  $\C \mathbb
P^n$. One can also view the next statement as a global analog of
the Nijenhuis-Woolf theorem \cite{NiWo,Si}.

\begin{theorem}\label{deform}
Let $J$ be a smooth almost complex structure in $\C^n$ satisfying Conditions~(\ref{size0}) and (\ref{size1}), where $\lambda$ is sufficiently small. Then there exists a positive constant $C$ such that for every $v \in \mathbb S$ there
is a $J$-complex line $L_J^{v} : \mathbb C \rightarrow \mathbb C^{n}$ satisfying~:
\begin{eqnarray}
\label{growth}
\sup_{\zeta \in \mathbb C}\parallel L_J^{v}(\zeta) - \zeta v \parallel \leq C \lambda.
\end{eqnarray}
Furthermore
$$
\mathbb C^n \backslash \mathbb B =\{L_J^v(\zeta),\ v \in \mathbb S,\ \zeta \in \mathbb C\}.
$$
\end{theorem}

\bigskip

Denote by $L_p(\C,\C^n)$ the standard space of Lebesgue
$p$-integrable maps from $\C$ to $\C^n$ and by $\mathcal
C^{m,\gamma}(\C,\C^n)$ the space of $m$ times continuously
differentiable maps from $\C$ to $\C^n$, with $\gamma$-Holderian
partial derivatives of order $m$ (globally on $\C$), both equipped
with the standard norms. Recall that

\begin{eqnarray}
\label{norm0}
\parallel z
\parallel_{\mathcal C^{m,\gamma}(\C,\C^n)} = \parallel z
\parallel_{\mathcal C^m(\C)} + \sum_{\vert \alpha \vert = m}
\sup_{\zeta,\zeta' \in \C} \frac{\parallel D^\alpha z(\zeta) -
D^\alpha z(\zeta')\parallel }{\vert \zeta - \zeta'\vert^\gamma}
\end{eqnarray} As usual, the norm on the intersection $L_p(\C,\C^n) \cap
\mathcal C^{m,\gamma}(\C,\C^n)$ is the sum of the norms on
$L_p(\C,\C^n)$ and $\mathcal C^{m,\gamma}(\C,\C^n)$. Recall also
that the Cauchy-Green transform on the complex
plane $\C$ is defined by

$$Tz(\zeta) = \frac{1}{2\pi i}\int\int_{\mathbb C} \frac{z(\tau) d\tau \wedge
  d\overline\tau}{\tau - \zeta}.$$

We need the following classical property of the Cauchy-Green integral (see for
instance \cite{Ve}, p.63):

\begin{proposition}
\label{reg}
For every $p \in [1,2[$, for every integer $m \geq 0$ and every real number $\gamma \in ]0,1[$
$$T:L_p(\C,\C^n) \cap \mathcal C^{m,\gamma}(\C,\C^n) \rightarrow L_p(\C,\C^n) \cap \mathcal C^{m+1,\gamma}(\C,\C^n)$$
is a linear bounded operator satisfying
$$(Th)_{\overline\zeta} = h.$$
\end{proposition}

\vskip 0,2cm \noindent{\sl Proof of Theorem~\ref{deform}}. This is
based on a bijective correspondence between certain classes of
$J$-complex lines and standard
complex lines. The crucial observation is the following. Let $L:
\C \rightarrow \C^n$, $L: \zeta \mapsto z(\zeta)$ be a
$J$-holomorphic map. Then the map $L$ is a solution of the
following quasi-linear elliptic system of partial differential
equations:

\begin{eqnarray}
\label{CR}
z_{\overline\zeta} - A(z) { \overline z}_{\overline\zeta} = 0
\end{eqnarray}
where $A(z)$ is a  complex $(n \times n)$-matrix function. It
corresponds to the matrix representation of the real endomorphism
$(J(z) + J_{st})^{-1}(J(z) - J_{st})$ of $\R^{2n}$; it is easy to
check that this endomorphism is anti-linear with respect to
$J_{st}$ (cf. for instance \cite{NiWo}). This equation is equivalent to

\begin{eqnarray}
\label{CR1} (z - T A(z) {\overline
z}_{\overline\zeta})_{\overline\zeta} = 0
\end{eqnarray}
that is to the $J_{st}$-holomorphicity of the map $w: \zeta
\mapsto z(\zeta) - T A(z(\zeta)) {\overline z}_{\overline
\zeta}(\zeta)$ on $\mathbb C$.

Consider the $J_{st}$-holomorphic map  $L^{v}: \zeta
\longrightarrow  \zeta v$, $\zeta \in \C$. Our goal is to
consider $J$-holomorphic maps defined on $\C$ close enough to
$L^{v}$. Since the map $L^{v}$ is not bounded and so does not
belong to the space $C^{1,\gamma}(\C,\C^n)$ equipped with the
standard norm introduced above, we need to modify slightly this norm
in order to enlarge this space. Fix $\gamma \in ]0,1[$ and
consider a continuous map $z:\C \longrightarrow \C^n$. We define
its weighted ${\mathcal C}^0$ norm by

\begin{eqnarray}
\label{sup} \parallel z \parallel_w:= \parallel z(\zeta)(1 + \vert
\zeta \vert^2)^{-1/2}\parallel_{{\mathcal C}^0(\C,\C^n)}.
\end{eqnarray}
For a function $z:\C \longrightarrow \C^n$ of class ${\mathcal
C}^1$ with the $\gamma$-H\"olderian first order partial
derivatives we define the weighted ${\mathcal C}^{1,\gamma}$-norm
by
\begin{eqnarray}
\label{norm1} \parallel z \parallel_{C^{1,\gamma}_w(\C,\C^n)}:=
\parallel z \parallel_w + \parallel z_\zeta
\parallel_{{\mathcal C}^{0,\gamma}(\C,\C^n)} + \parallel
z_{\overline\zeta}
\parallel_{{\mathcal C}^{0,\gamma}(\C,\C^n)}
\end{eqnarray}
Thus we just add the weight to the ${\mathcal C}^0$-term in the
standard norm (\ref{norm0}). We denote by ${\mathcal
C}^{1,\gamma}_w(\C,\C^n)$ the space of the above maps $z: \C
\longrightarrow \C^n $ equipped with this norm. Then obviously
\begin{itemize}
\item[(i)] The space ${\mathcal C}^{1,\gamma}_w(\C,\C^n)$ is
Banach \item[(ii)] For every $z \in {\mathcal
C}^{1,\gamma}_w(\C,\C^n)$ one has
\begin{eqnarray}
\label{norms} \parallel z
\parallel_{C^{1,\gamma}_w(\C,\C^n)} \leq
\parallel z \parallel_{C^{1,\gamma}(\C,\C^n)}.
\end{eqnarray}
\end{itemize}
In particular the identity map
$$id:{\mathcal C}^{1,\gamma}(\C,\C^n) \longrightarrow {\mathcal
C}^{1,\gamma}_w(\C,\C^n)$$ is continuous.

For $\varepsilon_0 > 0$  small enough consider the open subset $\mathcal U_{\varepsilon_0}$ of the space
${\mathcal C}^{1,\gamma}_w(\C,\C^n)$ defined by~:
$$
\mathcal U_{\varepsilon_0}:=\{z \in \mathcal
C^{1,\gamma}_w(\C,\C^{n}) / \|z - L^{v}\|_{\mathcal
C^{1,\gamma}_w(\C,\C^n)}< \varepsilon_0, \|v\|=1\}.
$$
We point out that $\mathcal U_{\varepsilon_0}$ is independent of $v$ which
runs over the unit sphere in the above definition.

Given  $\theta > 1$ from Condition~(\ref{size0}) fix $p \in ]1,2[$
such that $\theta p  > 2$. Then $A(z)\bar{z}_{\bar{\zeta}} \in L_p(\mathbb C,\mathbb C^n) \cap \mathcal C^{0,\gamma}(\mathbb C,\mathbb C^n)$ for every $z \in \mathcal U_{\varepsilon_0}$. Thus the operator~:
$$
\begin{array}{ccccc}
\Phi_J &:& \mathcal U_{\varepsilon_0} \subset \mathcal
C^{1,\gamma}_w(\C, \C^n) & \rightarrow & \mathcal C^{1,\gamma}_w(\C, \C^n)\\
 & & z &\mapsto & w = z - T A(z) {\overline z}_{\overline\zeta}.
\end{array}
$$
is correctly defined by Proposition \ref{reg} if $\varepsilon_0  > 0$ is small
enough. In what follows we always assume that this assumption is satisfied.

\begin{proposition}\label{oper}
If $\lambda$ given by Conditions (\ref{size0}), (\ref{size1}) is
sufficiently small then the operator $\Phi_J$ is a smooth local $\mathcal C^1$
diffeomorphism from $\mathcal U_{\varepsilon_0}$ to
$\Phi_J(\mathcal U_{\varepsilon_0})$. Moreover $L^{v} \in
\Phi_J(\mathcal U_{\varepsilon_0})$ for every $v \in \mathbb S$.
\end{proposition}

\noindent{\sl Proof of Proposition~\ref{oper}.} First we show   that $\Phi_J$ is a small $\mathcal C^1$ deformation of the identity on $\mathcal U_{\varepsilon_0}$.

\begin{lemma}\label{estim}
There exists a positive constant $D_1(\varepsilon_0)$ such that
for every
$z \in \mathcal U_{\varepsilon_0}$~:
\begin{eqnarray}\label{size6}
\parallel \Phi_J(z) - z\parallel_{\mathcal C^{1,\gamma}(\C,\C^n)} \leq D_1
\lambda.
\end{eqnarray}
\end{lemma}
We point out that we employ the standard non-weighted norm in
(\ref{size6}).

 \noindent{\sl Proof of Lemma~\ref{estim}}. Indeed, in view of
Condition~(\ref{size0}) and the choice of $p$ we have the estimate
\begin{eqnarray}
\label{size7}
\parallel A(z){\overline z}_{\overline \zeta} \parallel_{L_p(\C,\C^n)} + \parallel A(z){\overline z}_{\overline \zeta}
\parallel_{\mathcal C^{0,\gamma}(\C,\C^n)} \leq D_2\lambda
\end{eqnarray}
for any $z \in \mathcal U_{\varepsilon_0}$, for some $D_2>0$.
So Proposition \ref{reg} implies the statement of Lemma. \qed

\vskip 0,2cm
Next we  prove that the Fr\'echet
derivative $\dot \Phi_J$ of $\Phi_J$ is close to identity as an operator. The
derivative of $\Phi_J$ at $z \in \mathcal U_{\varepsilon_0}$ is defined by :
$$
\begin{array}{ccccc}
\dot \Phi_J(z) & : &  C^{1,\gamma}_w(\C,\C^n) & \rightarrow
& \mathcal C^{1,\gamma}_w(\C,\C^n)\\
 & & \dot z & \mapsto & \dot w=\dot z -TA(z) {\overline{\dot z}}_{\overline \zeta} - T(B(\dot z){\overline z}_{\overline \zeta})
 \end{array}
 $$
 where $B$ is given by
 $$
 B(\dot z) = [J(z)+J_{st}]^{-1}DJ(z)(\dot z) - [J(z)+J_{st}]^{-1}DJ(z)(\dot z)[J(z)+J_{st}]^{-1}[J(z)-J_{st}].
 $$

\begin{lemma}\label{der}
There exists a positive constant $D_3(\varepsilon_0)$ such that for every $z \in \mathcal U_{\varepsilon_0}$~:
\begin{equation}
\||\dot \Phi_J(z)-Id\|| \leq D_2 \lambda.
\end{equation}
\end{lemma}

\noindent{\sl Proof of Lemma~\ref{der}.} It follows from Condition
(\ref{size0})  that there exist positive constants $D_4, D_5$
such that for every $z \in \mathcal
U_{\varepsilon_0}$ and for every $\dot z \in \mathcal
C^{1,\gamma}_w(\C,\C^{n})$ we have~:

\begin{equation}\label{size8}
\parallel B(\dot z){\overline z}_{\overline \zeta} \parallel_{L_p(\C,\C^n)} +  \parallel B(\dot z){\overline z}_{\overline \zeta}
\parallel_{\mathcal C^{0,\gamma}(\C,\C^n)}
\leq D_4 \lambda  \|\dot z\|_{\mathcal C^{1,\gamma}_w(\C,\C^{n})},
\end{equation}
and

\begin{equation}\label{size9}
\parallel A(z){\overline {\dot z}}_{\overline\zeta} \parallel_{L_p(\C,\C^n)} +  \parallel A( z){\overline {\dot z}}_{\overline \zeta}
\parallel_{\mathcal C^{0,\gamma}(\C,\C^n)}
\leq D_5 \lambda \|\dot z\|_{\mathcal C^{1,\gamma}_w(\C,\C^{n})}.
\end{equation}

Now  Proposition \ref{reg} implies the desired statement. \qed

\begin{remark}
We point out that the estimates~(\ref{size7}), (\ref{size8}), (\ref{size9}) and the existence of the Fr\'echet derivative of $\Phi_J$ rely deeply on the definition of $\mathcal U_{\varepsilon_0}$. Indeed every map $z \in \mathcal U_{\varepsilon_0}$ has the same asymptotic behaviour at infinity as $|\zeta|$. This allows to use Conditions (\ref{size0}) and (\ref{size1}).
\end{remark}

Proposition~\ref{oper} is now a consequence of Lemmas \ref{estim}
and \ref{der} provided $\lambda$ is sufficiently small. We just
need to justify the last statement, i.e. the condition $L^{v}
\in  \Phi_J(\mathcal U_{\varepsilon_0})$ for every $v \in \mathbb S$.
We must show that the equation~:
\begin{eqnarray}
\label{fixed} z = T A(z){\overline z}_{\overline \zeta} + L^{v}
\end{eqnarray}
admits a solution $z \in {\mathcal U_{\varepsilon_0}}$. Consider
the operator $Q: z \mapsto T A(z){\overline z}_{\overline \zeta} +
L^{v}$ defined on the closure $\overline {\mathcal
U_{\varepsilon_0}}$. It follows by Lemma \ref{estim} and
the estimate (\ref{norms}) that for $\lambda$ small enough the image
$Q(\overline {\mathcal U_{\varepsilon_0}})$ is contained in
$\overline {\mathcal U_{\varepsilon_0}}$. Furthermore Lemma \ref{der}
implies that $Q$ is a contracting map. So by the Fixed Point Theorem the
equation (\ref{fixed}) has a unique solution in $\overline
{\mathcal U_{\varepsilon_0}}$. \qed

\vskip 0,1cm We can conclude now the proof of
Theorem~\ref{deform}. According to Proposition~\ref{oper} if
$\lambda$ is sufficiently small then
for every $v \in \mathbb S$  the map $L^{v}$ belongs to
$\Phi_J(\mathcal U_{\varepsilon_0})$. Then the map
$L_J^v:=\Phi_J^{-1}(L^{v})$ is $J$- holomorphic on $\C$,
satisfying the condition (\ref{growth}).
Finally consider the evaluation map
$$
\begin{array}{ccccc}
ev_J & : & \mathbb R^+ \times \mathbb S & \rightarrow & \mathbb C^n\\
   &   & (t,v)  & \mapsto & L_J^v(t)\;
\end{array}
$$
smoothly depending on $J$ as a parameter. The map $ev_J$ is a smooth small
deformation of the map $ev_{J_{st}} : (t,v) \mapsto tv$. Let the projection $\pi:
\mathbb C^n\backslash \{ 0 \}  \rightarrow 
\mathbb R^+ \times \mathbb S$ be defined by $\pi(z) = (\parallel
z \parallel^2, \parallel z \parallel^{-2}z)$ 
so that $(ev_{J_{st}} \circ \pi):\mathbb C^n\backslash \{ 0 \}  \rightarrow
\mathbb C^n\backslash \{ 0 \} $ 
is the identity map. Now we extend the restriction $(ev_{J} \circ \pi)\vert_{\mathbb S}$ smoothly on the unit ball $ \mathbb B$ and obtain a smooth map $Ev_{J}: \mathbb C^n \rightarrow \mathbb C^n$ coinciding with   $(ev_{J} \circ \pi)$ on $\mathbb C^n \backslash   \mathbb B$.  It follows by the condition 
(\ref{growth}) that we have $\parallel  Ev_J (z) \parallel \rightarrow \infty$ as $\parallel z \parallel \rightarrow \infty$ i.e. this map is proper and so  is surjective. \qed

\section{Examples and remarks}

{\bf 5.1.}
It is not difficult to
produce examples of structures satisfying the assumptions of Theorem
\ref{MainTheorem}. One may consider a smooth diffeomorphism $F: \R^{2n}
\rightarrow \R^{2n}$ such that $F$ converges ``fast enough'' to the identity map  as $\sum_j \vert z_j \vert^2$ tends to infinity, and set
$J= F_*(J_{st})$.

The integer $K(n)$ in Theorem \ref{MainTheorem} depends on results of
G.Kiremidjian and L.Lempert. In the case $n \geq 3$ it follows by \cite{Kir} that
one can take $K(n) = n+3$. In the case $n=2$ the assumptions on the initial
regularity are more involved (see Proposition 13.1 of \cite{Le}).

We point out that  Charles L. Epstein and Yong Ouyang~\cite{ep-yo}
obtained another version of Theorem \ref{MainTheorem}. Roughly speaking they
require the convergence of $J$ and its partial derivatives up to the third
order to $J_{st}$ at infinity, with ``the third degree plus $\varepsilon$''  polynomial decrease.
Their approach is based on the study of the asymptotic behavior of sectional
curvatures near a pole in a complete K\"ahler manifold. Thus our assumption on
the asymptotic behaviour of the partial derivatives of $J$ up to the third
order are weaker, but we need assumptions on higher order derivatives.

The approach of I.V.Zhuravlev \cite{Z} is based on an explicit solution of the $\overline\partial$-equation in $\C^n$ by means of a suitable integral representation.
He obtains an analogue of Theorem \ref{MainTheorem} assuming that the norm
$\parallel J - J_{st} \parallel_{L^\infty(\C^n)}$ is small enough and the
matrix function $J - J_{st}$ admits certain second order Sobolev derivatives
in $L^p(\C^n)$ for suitable $p > 1$. This result requires a quite low
regularity of $J$. However Sobolev's type  condition of $L^p$ integrablity on
$\C^n$ is somewhow restrictive: there are obvious examples of functions
polynomially decreasing on $\C^n$,  which are not in $L^p(\C^n)$. So our result is independent from the results of \cite{ep-yo,Z}.
It would be interesting to find a general statement which would contain the
results of \cite{ep-yo,Z} and the result of the present paper as special
cases.

A strong result was obtained recently by E.Chirka \cite{Ch} who proved an analogue of Theorem \ref{MainTheorem} in the case of $\C^2$,
under the assumption that  the norm $\parallel J - J_{st}
\parallel_{L^\infty(\C^2)}$ is small enough and that $J - J_{st}$ is of class
$L^2$ on $\C^2$. His method is based on the study of global foliations of $\C^2$
by pseudoholomorphic curves and strongly uses the dimension two. He also obtained a statement on the existence of $J$-complex lines 
in the spirit of Theorem \ref{deform} under assumptions different from ours; his result can not be applied in our situation.
We believe that the  approaches of the present work and the works \cite{Ch,ep-yo,Z} will be useful for a further work concerning the  natural  questions on deformations of Stein structures.

{\bf 5.2.} We point out that there exist integrable almost complex structures on
$\R^{2n}$ without non-constant plurisubharmonic functions (see \cite{DiSi}). The simplest
examples come from the complex projective space $\C \mathbb P^n$.

\begin{proposition}\label{example}
There exists an integrable almost complex structure $J$ on $\R^{2n}$ such that every $J$-plurisubharmonic function $u: \R^{2n} \rightarrow [-\infty,+\infty[$ is constant.
\end{proposition}
\proof Consider the complex projective space $\C \mathbb P^n$ and the affine space
$\C^n$ equipped with the standard complex structures $J_{st}^\mathbb P$ and $J_{st}$
respectively. Let also $\pi: \C \mathbb P^n \rightarrow \C^n$ be the canonical
projection. Then $\C \mathbb P^n = \pi^{-1}(\C^n) \cup S$ where $S$ ("the pullback of
infinity") is a smooth compact complex hypersurface in $\C \mathbb P^n$. After an arbitrary small perturbation of $S$, given by a global diffeomorphism $\Phi$ of $\mathbb C \mathbb P^n$ we obtain a smooth compact
real submanifold $\tilde S:=\Phi(S)$ of codimension 2 in $\C \mathbb P^n$.

\vskip 0,5cm
For such a
general perturbation, the manifold $\tilde S$ will be generic almost
everywhere, meaning that the complex linear span of its tangent space in almost
every point is equal to the tangent space of $\C \mathbb P^n$ at this point. According
to well known results (see for instance \cite{Ka}) such a manifold is
a removable singularity for any plurisubharmonic function. More precisely, if
$u$ is a plurisubharmonic function on $\C \mathbb P^n \backslash \tilde S$ then there
exists a plurisubharmonic function $\tilde u$ on $\C \mathbb P^n$ such that
$\tilde u \vert_{\C \mathbb P^n \backslash \tilde S} = u$. Therefore any function
which is plurisubharmonic (with respect to $J_{st}^{\mathbb P}$) on $\C \mathbb P^n \backslash
\tilde S$ is constant.
Then the complex structure $J:= (\pi \circ \Phi^{-1})_*(J_{st}^\mathbb P)$ satisfies the statement of Proposition~\ref{example}.

\vskip 0,2cm
\noindent{\it Acknoweldgements.} We thank E.M.Chirka for bringing our attention to the subject of the present paper and to the  works
\cite{Ch,Z}. We are grateful to Ch.Epstein for bringing our attention to the  article~\cite{ep-yo}.

\end{document}